# HYPERGRAPH VERSALS

*Vance Faber and Noah Streib*

*revision: December 6, 2015*

**Introduction**. Let $H(V,E)$ be a hypergraph on $n$ vertices $V$ and $m$ edges $E$ with the property that no edge contains another. The rank $r$ of $H$ is the maximum cardinality of an edge. We say a set of vertices $S$ is a *versal* for the edge $e$ if for every other edge $f$

$$|e| + |S \cap e| < |f| + |S \cap f|.$$

Let $L(e) = L(H,e)$ be the set of all versals for the edge $e$. Let $Z(H) = \cup \{L(e) \mid e \in E\}$ be the set of all versals. We shall show the following.

**Main Theorem.** The cardinality of $Z(H)$ is at least $n+1$ unless $H$ is

a) the set of all singletons $S_n$,

b) the set of all complements of singletons $\widetilde{S}_n$,

c) the graph $C_4$ (the four-cycle.)

In those cases, the cardinality is $n$.

**Motivation**. The question is (weakly) connected to certain parallel algorithms [2] where it is necessary to have a lower bound on the probability that a function from the vertices of a hypergraph into a set of integers will have a single edge where the sum is minimum. The idea of a versal corresponds to a function with a unique minimum which takes the value 2 on the versal and 1 off it. The question in the function formulation was given as an open problem in [1].

**Proof in the uniform case**. To simplify the arguments, we start with a uniform hypergraph. A uniform hypergraph is one with all edges of identical cardinality $r$. In a uniform hypergaph, the condition for a versal is simplified to

$$|S \cap e| < |S \cap f|$$

for every edge $f \neq e$.

*Definition.* Given a hypergraph $H$ we denote by $\widetilde{H}$ the *reflection* hypergraph whose edges are the complements of the edges of $H$. We shall use $\bar{e}$ to denote the complement of the set $e$.

*Lemma 1.* If $H$ is a uniform hypergraph the cardinalities of $Z(H)$ and $Z(\widetilde{H})$ are identical.

*Proof.* There is a one-to-one map from $Z(H)$ to $Z(\widetilde{H})$. Suppose $S$ is a versal for the edge $e$ and $f$ is some other edge. Then since $|e| = r = |f|$

$$|\bar{S} \cap \bar{e}| = |V \setminus (S \cup e)| = n - |S \cup e| = n - |S| - |e| + |S \cap e|$$

$$< n - |S| - |f| + |S \cap f| = |\bar{S} \cap \bar{f}|$$

so $\bar{S}$ is a versal for $\bar{e}$.

Consequently, to prove the Main Theorem in the uniform case, we need only focus on the case $2r \leq n$.

*Definition.* Fix $r$. An $m$-star $T_m$ (Figure 1 shows $T_4$) is a uniform hypergraph with $m > 1$ edges and rank $r$ such that all edges share the same $r - 1$ vertices called the *core* of the star. (Stars are a subclass of sunflowers [3].) An $m$-binary star $B_m$ (Figure 2 shows $B_4$) is the union two $m$-stars such that the intersection of the cores has cardinality $r - 2$ and each vertex outside the core of one is identical to a distinct vertex outside the core of the other. Note that $B_2 = C_4$, the 4-cycle. There are two vertices that are not shared between the two stars. They form the symmetric difference of the cores.

*Definition.* A versal $S$ for the edge $e$ is called *null* if $S \subseteq \bar{e}$. Let $Z'(H)$ be the set of all null versals. Note that the condition for $S$ to be a null versal can be rewritten as $S \subset \bar{e}$ and for every edge $f \neq e$, $S \cap f \neq \emptyset$. In particular, $\bar{e}$ is a null versal for $e$.

*Theorem 2.* If $H$ is uniform with rank $r$, $2r \leq n$, $m > 1$ then the cardinality of $Z'(H)$ is at least $n + 1$ unless $H$ is

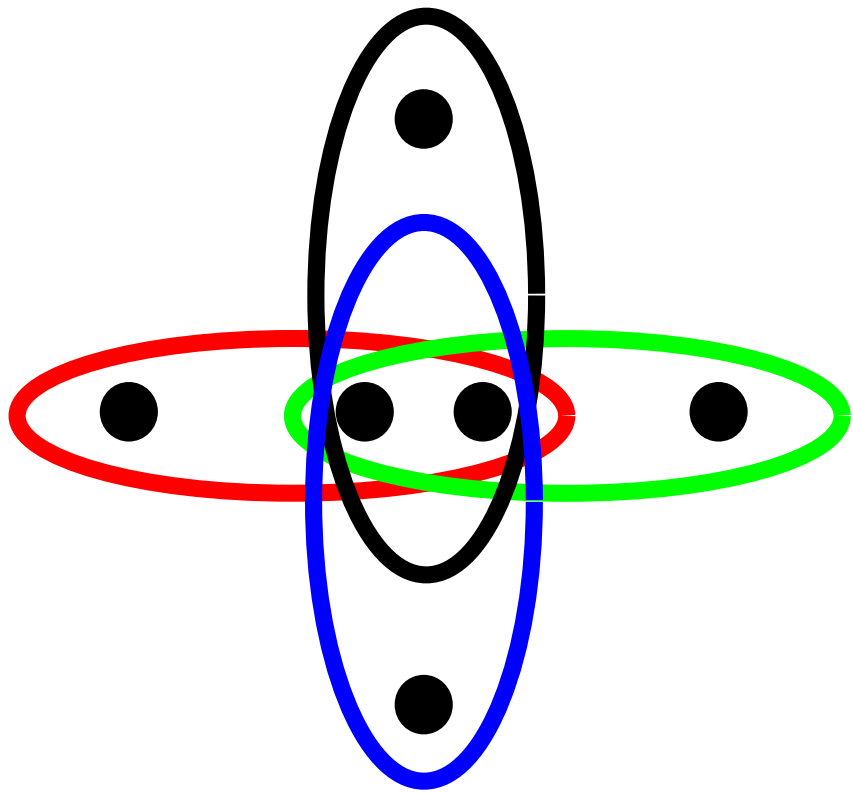

Figure 1. 4-star with r=3.

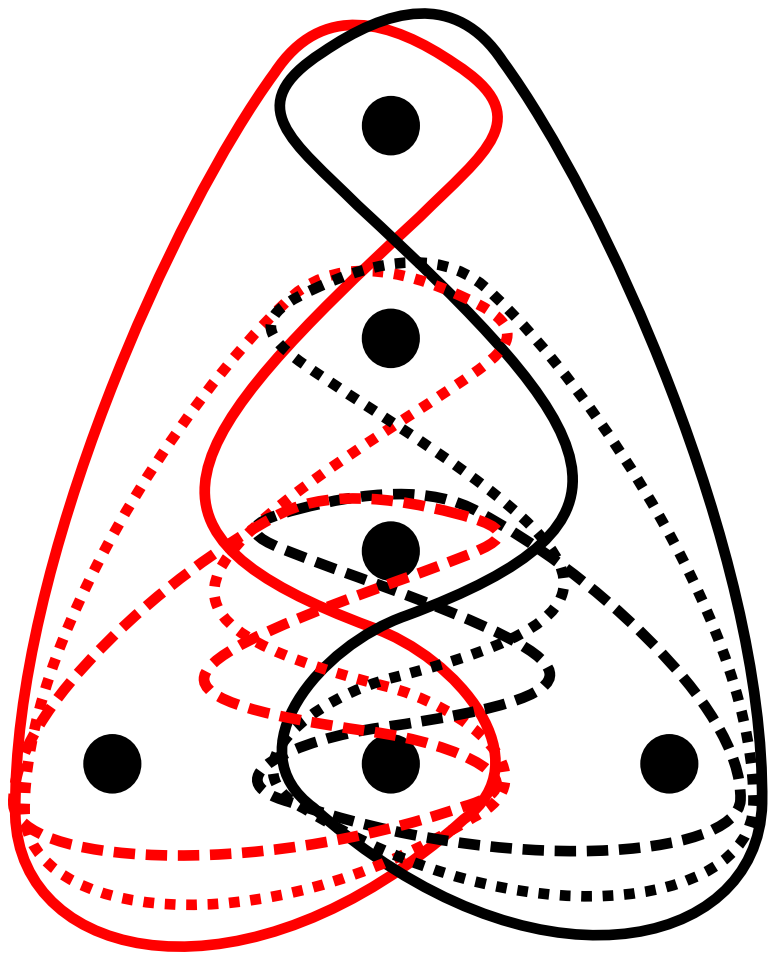

Figure 2. 4-binary star with r=3.

a) the singleton hypergraph $S_n$,

b) the star $T_m$,

c) the binary star $B_m$.

In case b) the cardinality is $m = n - r + 1$ and in the other cases, the cardinality is $n$.

Before proceeding with the proof, we need some definitions and two technical lemmas.

*Definition.* Let $S$ be a versal for the edge $e$. We say a vertex $v \notin S \cup e$ is *free* with respect to the edge $e$ and the versal $S$. Note that $S \cup \{v\}$ is also a versal for $e$ (since $(S \cup \{v\}) \cap e = S \cap e$) and if there are $p$ free vertices with respect to $e$ and $S$ then $|L(e)| \geq 2^p$.

*Lemma 3.* Let $e$ be an edge in a (not necessarily uniform) hypergraph $H$ of smallest cardinality $\rho$. Suppose that $v$ is a vertex in $\bar{e}$ that is not free with respect to $e$. Then there must exist an edge $f_v$ in $H$ of cardinality $\rho$ such that $v$ is its only vertex in $\bar{e}$.

*Proof.* Note that if every edge of cardinality $\rho$ that contains $v$ also contains another vertex in $\bar{e}$ then $v$ is free with respect to $e$ and $\bar{e} \setminus \{v\}$.

*Lemma 4.* Suppose $H$ has edges $e$ and $f$ and an edge $g$ of smallest cardinality $\rho$ in $H$ and a vertex $u \in f$ with

a) degree of $u$ is 1

b) $f \setminus \{u\} \subseteq e$.

If $|f \cap g| < \rho - 1$ then there exists a free vertex with respect to $g$.

*Proof.* First, we write the disjoint union

$$f = [f \setminus (e \cup g)] \dot\cup [(f \cap e) \setminus g] \dot\cup [(f \cap g) \setminus e] \dot\cup [f \cap g \cap e].$$

Since degree of $u$ is 1, then $u$ is in $f$ but not $e \cup g$. Since in addition $f \setminus \{u\} \subseteq e$ then $f \setminus (e \cup g) = \{u\}$. If $(f \cap e) \setminus g$ were empty then

$$f = \{u\} \, \dot\cup \, [(f \cap g) \setminus e] \, \dot\cup \, [f \cap g \cap e)] = \{u\} \, \dot\cup \, (f \cap g).$$

But then $\rho \leq |f| = 1 + |f \cap g|$ and thus $|f \cap g| \geq \rho - 1$, a contradiction, so $(e \cap f) \setminus g$ is not empty. Then $u$ is free with respect to $g$ and the versal $\bar{g} \setminus \{u\}$.

To make the proof of Theorem 2 easier to understand, we need to encapsulate a few more concepts.

*Definition.* Let $n \geq 2r$. A uniform hypergraph $F(e)$ with $m = n - r + 1$ edges is a *flag* with *pole* $e$ and *pennants* $e_1, e_2, \ldots, e_{n-r}$ if $|e_j \cap e| = r - 1$ for all $j$ and if $v$ is the unique vertex in $e_j - e$ then $d(v) = 1$.

*Lemma 5.* The following statements about a hypergraph $F$ with $m = n - r + 1$ are equivalent:

a) $F$ has more than one pole;

b) $F$ is a star;

c) Every edge is a pole.

*Proof.* It is sufficient to prove that a) implies b). Let $e$ and $f$ be poles of flags. Since $m = n - r + 1$, each must be a pennant of the other. The remaining edges must be pennants of both edges and so must contain the set of $r - 2$ vertices common to both $e$ and $f$. Thus $F$ is a star.

*Lemma 6.* Let $H$ be an $r$-uniform hypergraph on $n$ vertices and $m$ edges with $n \geq 2r$ and $n - r + 1 \geq m \leq n$. Then every edge of $H$ is the pole of a flag if and only if $H$ is an $(n - r + 1)$-star or $m = n = 2r$ and $H$ is the $m$-binary star.

*Proof:* Let $F$ be a flag in $H$ with pennants $P$ and pole $z$. Disregard the remaining edges in $H$ for the moment---the proof proceeds by examining the ways in which these other edges might exist. Let $N(e)$ be the number of pennants missing from an edge $e$ in order to make it the pole of a flag. For example, $N(z) = 0$ for the pole of $F$. If every edge is a pole, $N(e) = 0$ for all edges $e$. We shall use 3 properties of $F$.

*Prop 1.* For each $p$ in $P$, $N(p) = n - r - k$ and there are $k - 1$ other pennants in $P$ that map to the same $(r - 1)$-set of $F$.

This is because the edges $p_1$ and $p_2$ in $P$ are pennants of each other if and only if they map to the same $(r-1)$-set of $F$, and each edge requires $n-r$ pennants to be a pole of a flag.

*Prop 2.* Let $e$, $f$ be edges. If the intersection of $e$ and $f$ is exactly $r-2$ then any edge that is a pennant of each must share that $(r-2)$-set and contain one of the two vertices in $e \setminus f$ and one of the two in $f \setminus e$. If the intersection is strictly smaller than $r-2$, then there is no edge that is a pennant of each.

*Prop 3.* Let $e$, $f$ in $P$ with intersection $r-2$ and $N(e), N(f) > 0$. Then, starting with the graph that only contains the edges of $F$, there is a unique edge that reduces both $N(e)$ and $N(f)$.

This is because the edge needs to be a pennant of each, and by Prop 2, there are only four such edges. One is $z$, which already exists as a pennant of each. Two others intersect $z$ in an $(r-1)$-set, and hence only reduce one of $N(e)$ or $N(f)$. The remaining edge, which contains the vertex in $e \setminus z$ and the vertex in $f \setminus z$ reduces both.

Now the proof breaks into three cases.

Case 1. There is a $p$ in $P$ with $2 \leq N(p) \leq n-r-2$. We shall show this leads to a contradiction.

Pick such a $p$. There are $2 \leq k \leq n-r-2$ other $p$ in $P$ that map to the same $(r-1)$-set of $F$, where $N(p) = n-r-k$ Call this set $P'$. Thus, there are $n-r-k$ pennants in $P \setminus P' = P''$. Notice that $N(p) \geq n-r-k$ for all $p$ in $P''$ since no edge in $P'$ is a pennant for any edge in $P''$.

Let $p'$ in $P'$ and $p''$ in $P''$. We need to add $m - (n-r+1)$ more edges to our graph. By Prop 3, there is a unique edge $e$ that is a pennant for $p'$ and $p''$ and reduces $N(p')$ and $N(p'')$. Add e to the graph. Now $N(p') = k-1 > 0$, $N(p'') \geq n-r-k-1 > 0$, and $N(p') + N(p'') \geq n-r-2 \geq r-2$. Since there are $m-n+r-2 \leq r-2$ remaining edges to add and no more edges can decrease both $N(p')$ and $N(p'')$, every remaining edge must be a pennant of one or the other. Notice that if we did not select $e$, then we would need to add at least $r-1$ edges without using an edge to that reduced both $N(p')$ and $N(p'')$, which makes it impossible to make each a pole.

Notice that $|P'| \geq 2$ and $|P''| \geq 2$. Therefore there exist $q'$ in $P' \setminus \{p'\}$ and $q''$ in $P'' \setminus \{p''\}$. Applying the same argument, we must add $f$ to our graph, the unique edge that reduces both $N(q')$ and $N(q'')$. However, $f$ is neither a pennant for $p'$ or $p''$, implying that we can't make poles of flags out of each of $p'$, $p''$, $q'$, and $q''$.

Case 2. $N(p) = 1$ for $n-r-1$ elements of $P$ and $N(q) = n-r-1$ for the other element of $P$. In this case, we shall show that $H$ is a $2r$-binary star.

We must add $m-(n-r+1) \leq r-1$ edges to the graph we have so far. Thus, in order to make $q$ a pole of a flag, we need $n = 2r$ and each remaining edge a pennant of $q$. This further implies that each remaining edge reduces both $N(q)$ and $N(p)$ for distinct elements $p$ of $P \setminus q$. By Prop 3, these edges are uniquely determined. As it can be checked, adding these $r-1$ edges creates a flag decomposition. Further, the elements of $P \setminus q$, along with $z$, induce an $r$-star, and $q$ along with the new edges also induce an $r$-star, and their union is an $2r$-binary star.

Case 3. $N(p) = 0$ for all $p$ in $P$. That is, $F$ is an $(n-r+1)$-star. We shall show that $H = F$.

Let $R$ be the core of $F$. Let $e$ be an edge not in $F$. If the intersection of $e$ and $R$ is strictly less than $r-2$, then none of the edges in $F$ are pennants for $e$. Since there are only $m-(n-r+1)-1 \leq r-2$ edges remaining, and $e$ requires $n-r \geq r$ more pennants, we can't make $e$ into a pole. So the intersection of $e$ and $R$ must be exactly $r-2$ (since all edges that intersect $R$ in a set of size $r-1$ are in $F$ already). Then $e$ gains exactly two pennants from edges in $F$. However, these two pennants use the same vertex in the complement of $e$, meaning that $e$ still needs $r-1$ more pennants, which we can't make with our remaining $r-2$ edges. Therefore, if $H$ has a flag decomposition, it cannot contain any edges not in $F$. Since an $(n-r+1)$-star is itself a flag decomposition, it is only such decomposition obtained from this case.

Now we can prove Theorem 2.

*Proof of Theorem 2.* We have $n \geq 2r$. We break the proof into these cases and subcases:

Case I. $m > n$.

Case II. $m \leq n-r$.

   Case II.a. $n > 2r$.

   Case II.b. $n = 2r$.

      Case II.b.1. $m < r$.

Case II.b.2. $r \leq m \leq n - r = 2r$.

Case III. $m = n - r + k$ with $1 \leq k \leq r$.

   Case III.a. $n > 2r$.

   Case III.b. $n = 2r$

Now for the proof.

Case I. $m > n$. As we remarked above, for each edge $e$, $\bar{e}$ is a null versal so this case is trivial.

Case II.a. $m \leq n - r$, $n > 2r$. Given $e$ there exists a subset $B$ of $\bar{e}$ of cardinality at most $m - 1$ containing a vertex in each edge other than $e$. Then every subset $S$ of $\bar{e} \setminus B$ is a null versal for $e$. Thus

$$|Z'(H)| \geq m 2^{n-m-r+1}.$$

It follows that $|Z'(H)| \geq 2(n - r) \geq n + 1$ since $n \geq 2r + 1$.

Case II.b.1. $m \leq n - r$, $n = 2r$, $r \leq m \leq 2r$. Since $|\bar{e}| = r \geq m + 1 > m - 1$, we can find a subset $B$ of $\bar{e}$ of cardinality at most $m - 1$ containing a vertex in each edge other than $e$. Then $\bar{e} \setminus B$ contains at least $r - (m - 1) = r - m + 1$ vertices free with respect to $e$ and $B$. Thus

$$|Z'(H)| \geq m 2^{r-m+1}.$$

It follows that $|Z'(H)| \geq 4(r - 1) \geq n + 1$ since $n = 2r$ and $r \geq m + 1 \geq 3$.

Case II.b.2. $m \leq n - r$, $n = 2r$, $m < r$. If $H$ is a graph with 4 vertices and 2 edges then either the graph is a path or two disjoint edges neither of which have less than 5 versals. So we can assume $r \geq 3$. Since $|\bar{e}| = n - r = r$ while $|H| \setminus \{e\} = r - 1$ there is at least one free vertex $v$ with respect to $e$ in $\bar{e}$. This means that we have accounted for $n$ versals and we just need to find one more. Then $S = \bar{e} \setminus \{v\}$ will contain a second free element unless each of the $r - 1$ vertices in $S$ is a member of exactly one of the $r - 1$ edges in $|H| \setminus \{e\}$. Furthermore, no two of the edges in $|H| \setminus \{e\}$ can share $v$ or we can remove their vertices from $S$ and add back $v$ to obtain one additional versal. The vertex $v$ cannot be isolated or else the hypergraph $H(V \setminus \{v\}, E)$ has rank $r \geq 3$, number of vertices $n - 1 > r + 1$ and so must have at least $n$ versals. Each of those is a versal for $H$ and so adding on the vertex $v$ to any one gives at least $n + 1$ versals for $H$. This leads to the following configuration. Each of the edges $f \neq e$ has exactly one vertex in $\bar{e}$ except one edge $g$ which has the vertex $v$ and one other vertex $u$ which is not in any other edge. Note that we can construct three versals with respect to $e$ inside $\bar{e}$: $S, (S \cup \{v\}) \setminus \{u\}, S \cup \{v\}$. This gives the additional required versal and proves the lemma for this case.

Case III. $m = n - r + k$ with $1 \leq k \leq r$, $n \geq 2r$. For an edge $e$, let $P(e)$ be the set of pennants of the edge e. Each vertex in the complement of $e$ that does not have an associated pennant for $e$ is free with respect to $e$. Let $N(e)$ be the number of pennants needed to make $e$ into a pole. Now we assume by way of contradiction that $|Z'(H)| \leq n$. Let $q$ be the number of edge/vertex pairs such that the vertex is free for the edge. Then $m + q = n - r + k + q \leq |Z'(H)| \leq n$ so $q \leq r - k$. Furthermore, the number of edges that are poles of flags is at least $m - q$ and $m - q \geq (n - 2r) + 2k$ so there are least $(n - 2r) + 2k$ flags. At this point, we have shown that we have at least $2$ flags. However, by Lemma 5 if $k = 1$ and $H$ has at least 2 flags then $H$ is a star so we can assume that $k > 1$. Thus we have shown that there are at least 4 flags. We show that there are two flags with poles that are pennants of the each other. Let $e$ and $g$ be poles of flags. Each pole has $n - r$ pennants, and hence $2(n - r)$ between them (including ones counted twice). The intersection of their two sets of pennants must be as large as

$$P(e) \cap P(g) \geq 2(n - r) - m = n - r - k.$$

If $n > 2r$ then $P(e) \cap P(g) > r - k$. However, the number of edges in $H$ that are not poles is at most $q \leq r - k$. So there must be a pole $f$ that is in $P(e) \cap P(g)$. Then $|f \cap e| = r - 1$ and so $e$ is also a pennant of $f$. If $n = 2r$ and $k = r$ then since $m = 2r$ and there are $2k = 2r$ flags, every edge is the pole of a flag. By Lemma 6, $H$ is a binary star $B_{2r}$. Thus we can assume that $2 \leq k < r$ and $P(e) \cap P(g)$ is not empty. It follows that $|e \cap g| \geq r - 2$. If this intersection has cardinality $r - 1$, then $e$ and $g$ are pennants of each other. Otherwise, $|e \cap f| = r - 2$. Since $e$ and $g$ are poles, each vertex $v$ not in $e \cup g$ has two edges, one to be a pennant of $e$ and one to be a pennant of $g$. These two edges must be distinct or else the $r - 1$ vertices other than $v$ would be common to both $e$ and $g$. There are $n-(r+2)=r-2$ such vertices, meaning $2r-4$ such edges. Adding $e$ and $g$, we have $2r-2$ edges. The entire graph has $m=r+k\leq 2r-1$ edges. However, we require at least two more edges to make both $e$ and $g$ poles (they need to use the vertices in $e \setminus f$ and $f \setminus e$, respectively), a contradiction.

Thus we have flag poles $e$ and $f$ which are pennants of each other. Let $|P(e) \cap P(f)| = n-r-k+j$ with $j \geq 0$. The number of vertices outside of $e \cup f$ is $n-(r+1)$. Each of those vertices is either in a pennant shared by both $e$ and $f$ or it is in a different pennant for each of them. The cardinality of the set $Q$ of vertices not in pennants shared by both is

$$n-(r+1)-(n-r-k+j) = k-1-j.$$

so $j \leq k-1$. Suppose $j = k-1$ (recall that we are assuming $k > 1$.) In this case $Q$ is empty and the $n-r+1$ edges that we have counted for so far induce a star. Let $R$ be the core of the star. We must add $k-1$ additional edges to the star to form $H$. Let $g$ be the first such edge. Notice that $|g \cap R| \leq r-2$, as all edges that contain all of $R$ are in the star already. If $|g \cap R| < r-2$, then none of the edges in the star are pennants for $g$, and if $|g \cap R| = r-2$, then two of the edges in the star are pennants for $g$, but they use the same vertex in $\bar{g}$ (the other vertex in $R$). So in either case, $N(e) \geq n-r-1$, and there are just $k-2$ edges remaining to add to the graph. This means that the number of free vertices for $g$ will be at least

$$(n-r-1)-(k-2) = n-r-k+1 \geq r-k+1.$$

This contradicts the fact that the number of free vertex/edge pairs is at most $r-k$ and thus $j < k-1$.

Case III.a. $m = n-r+k$ with $2 \leq k \leq r$, $|P(e) \cap P(f)| = n-r-k+j$, $0 \leq j < k-1$, $n > 2r$. Edges $e$ and $f$ are pennants of each other. In this case $Q$ is non-empty. We have now tabulated the following number of edges

$$2+2(k-1-j)+n-r-k+j = n-r+k-j = m-j$$

so there are $j$ edges outside the union of the two flags $F(e)$ and $F(f)$. Notice that none of the pennants of $F(e)$ or $F(f)$ that contain a vertex in $Q$ can be pennants of any edge in $P(e) \cap P(f)$, and we have not identified enough edges to make an edge in $P(e) \cap P(f)$ into the pole of a flag. In order to decrease $N(g)$ for some $g$ in $P(e) \cap P(f)$, we need a pennant of $g$ that uses a vertex of $Q$. Such an edge can only decrease $N(g)$ for one $g$ in $P(e) \cap P(f)$. Since we only have $j$ remaining edges to add to our graph, we can make at most $j$ of the edges in $P(e) \cap P(f)$ into poles. This means our total number of flags is at most

$$m - P(e) \cap P(f) + j = n - r + k - (n - r - k + j) + j = 2k.$$

In this case where $n > 2r$, we require strictly greater than $2k$ poles so we have a contradiction.

Case III.b. $m = n - r + k$ with $2 \leq k \leq r$, $|P(e) \cap P(f)| = n - r - k + j$, $0 \leq j < k - 1$, $n = 2r$. The proof of Case III.a relied on the fact that we had strictly more than $2k$ poles which is not true in this case. Just before adding the final $j$ edges, we noticed that none of the edges in $P(e) \cap P(f)$ were yet poles of flags. Using $n = 2r$, we have $|P(e) \cap P(f)| = r - k + j$. Let $g$ be an edge in $P(e) \cap P(f)$. If we were to make $g$ into a pole with our final $j$ edges, notice that it would require adding $j$ edges, each of which included the vertex in $g$ that is not in $e \cap f$ and a distinct vertex in $Q$ (recall $|Q| = j$.) So in fact, we can only make one edge in $P(e) \cap P(f)$ into a pole. Since we have at most $r - k$ non-poles, this implies $j = 1$, and we add each of the edges we just described to the graph, which in this case is just one edge, which we call $g$. But note that $g$ itself must be a pole, since we have used all of our non-poles in $P(e) \cap P(f)$. By similar arguments, we find $|P(e) \cap P(f)|$ must be 1. However, then $|P(e) \cap P(f)| = r - k + j = r - k + 1 = 1$, which implies $r - k = 0$, a contradiction in this case.

*Theorem 7. If $H$ is uniform with rank $r$ then the cardinality of $Z(H)$ is at least $n + 1$ unless $H$ is either the graph $C_4$, the set of all singletons or the set of all complements of singletons. In those cases, the cardinality is $n$.*

*Proof.* The reader can verify that $Z(H)$ has cardinality exactly $n$ for the hypergraphs mentioned in the hypothesis. By Lemma 1 and Theorem 2, we need only show that the number of versals in a star or binary star is at least $n+1$. A binary star has $n$ null versals. If $r = 2$ the binary star is the same as $C_4$. If $r > 2$ then there is at least 1 vertex in the shared cores which can be added to any null versal to make additional versals. A star has $m = n - r + 1$ edges and $m$ null versals which are the complements of the edges. We get a versal from a null versal by including any subset of the common $r-1$ vertices. Thus there are

$$|Z(H)| = 2^{r-1}(n-r+1)$$

versals for the star. Since if $r \geq 3$ then $r \leq 2^{r-1} - 1$ we have

$$\frac{r}{2^{r-1}-1} \leq 1$$

$$r + \frac{r}{2^{r-1}-1} - 1 \leq r \leq n$$

$$r2^{r-1} - 2^{r-1} + 1 = (r-1)(2^{r-1}-1) + r \leq n(2^{r-1}-1) = n2^{r-1} - n$$

and finally

$$n+1 \leq n2^{r-1} + 2^{r-1} - r2^{r-1} = 2^{r-1}(n-r+1) = |Z(H)|.$$

**The main theorem.** Now we prove the Main Theorem.

*Proof of Main Theorem.* Suppose the theorem is false and let $H$ be a counterexample. Let $\rho$ be the cardinality of the smallest edge. Let $H'$ be the hypergraph consisting of only those edges of cardinality $\rho$. Note that $H'$ cannot be $C_4$, the set of all singletons or the set of all complements of singletons because for each of those hypergraphs any edge of larger cardinality will contain an edge already in the hypergraph. Suppose for the moment that $H'$ is not a star or a binary star.

Suppose first that $n \geq 2\rho$. Then by Theorem 2 there are at least $n+1$ null versals. Now we show that every null versal in $H'$ is a versal for $H$. Let $S$ be a null versal for the edge $e$ in $H'$ and let $f$ be an edge in $H$ with $e < |f|$. Then since $S$ is a subset of $\bar{e}$

$$|e| + |S \cap e| = |e| < |f| \leq |f| + |S \cap f|$$

so $S$ is a versal for $e$ in $H$. Thus $|Z(H)| \geq |Z(H')| \geq n+1$.

Now suppose that $n < 2\rho$. Then by Lemma 1, $\widetilde{H}'$ has at least $n+1$ null versals and their complements are versals for $H'$. Furthermore if $S$ is one such versal for an edge $e$ in $H'$, then because $\overline{S}$ is a null versal for $\overline{e}$, $\overline{S} \subseteq \overline{\overline{e}} = e$ so $S \supseteq \overline{e}$. In particular $S \cup e$ is the entire vertex set. Now let $f$ be an edge in $H$ with $|e| < |f|$. We can calculate

$$|e| + |S \cap e| = |e| + (e + |S| - |S \cup e|) = 2|e| + |S| - n < 2|f| + |S| - n$$

$$\leq 2|f| + |S| - |S \cup f| = |f| + |S \cap f|.$$

Thus $S$ is also a versal for $e$ in $H$.

Now suppose that $H'$ is a star or a binary star. Then any versal $S$ for $e$ in $H'$ contains $\overline{e}$. The proof in the previous paragraph shows that $S$ is also a versal for $e$ in $H$. This completes the proof of the theorem.

**References**.

1. Vance Faber, David Harris, "Tight bounds for the isolation lemma and min-unique functions" submitted.

2. https://en.wikipedia.org/wiki/Isolation_lemma

3. https://en.wikipedia.org/wiki/Sunflower_(mathematics)